\documentclass[12pt]{article}

\usepackage{amsmath, epsfig, cite}
\usepackage{amssymb}
\usepackage{amsfonts}
\usepackage{latexsym}
\usepackage{float}
\usepackage{color}

\newtheorem{thm}{Theorem}[section]

\newtheorem{lem}[thm]{Lemma}

\numberwithin{equation}{section}

\newcommand{\qed}{{\hfill$\square$}\medskip}

\begin{document}

\begin{center}
{\large\bf A $p$-adic analogue of Chan and Verrill's formula for $1/\pi$}
\end{center}

\vskip 2mm \centerline{Ji-Cai Liu}
\begin{center}
{\footnotesize Department of Mathematics, Wenzhou University, Wenzhou 325035, PR China\\
{\tt jcliu2016@gmail.com } \\[10pt]
}
\end{center}


\vskip 0.7cm \noindent{\bf Abstract.}
We prove three supercongruences for sums of Almkvist--Zudilin numbers, which confirm some conjectures of Zudilin and Z.-H. Sun. A typical example is the Ramanujan-type supercongruence:
\begin{align*}
\sum_{k=0}^{p-1}\frac{4k+1}{81^k}\gamma_k \equiv \left(\frac{-3}{p}\right)p\pmod{p^3},
\end{align*}
which is corresponding to Chan and Verrill's formula for $1/\pi$:
\begin{align*}
\sum_{k=0}^{\infty}\frac{4k+1}{81^k}\gamma_k=\frac{3\sqrt{3}}{2\pi}.
\end{align*}
Here $\gamma_n$ are the Almkvist--Zudilin numbers.


\vskip 2mm
\noindent{\it MR Subject Classifications}: 11A07, 11B65, 11Y55, 05A19

\section{Introduction}
For $n\ge 0$, the following sequence:
\begin{align*}
\gamma_n=\sum_{j=0}^{n}(-1)^{n-j}\frac{3^{n-3j}(3j)!}{(j!)^3}{n\choose 3j}{n+j\choose j}
\end{align*}
are known as Almkvist--Zudilin numbers (see \cite{az-b-2007} and A125143 in \cite{sloane-2020}).
This sequence appears to be first recorded by Zagier \cite{zagier-b-2009} as integral solutions to Ap\'ery-like recurrence equations.

These numbers also appear as coefficients of modular forms.
Let $q=e^{2\pi i \tau}$ and
\begin{align*}
\eta(\tau)=q^{1/24}\prod_{n=1}^{\infty}(1-q^n)
\end{align*}
be the Dedekind eta function. Chan and Verrill \cite{cv-mrl-2009} showed that if
\begin{align*}
t_3(\tau)=\left(\frac{\eta(3\tau)\eta(6\tau)}{\eta(\tau)\eta(2\tau)}\right)^4\quad
\text{and}\quad F_3(\tau)=\frac{\left(\eta(\tau)\eta(2\tau)\right)^3}{\eta(3\tau)\eta(6\tau)},
\end{align*}
and $|t_3(\tau)|$ is sufficiently small, then
\begin{align*}
F_3(\tau)=\sum_{n=0}^{\infty}\gamma_n t_3^n(\tau).
\end{align*}
They also constructed some new series for $1/\pi$ in terms of the numbers $\gamma_n$, one of the typical examples is the following formula \cite[Theorem 3.14]{cv-mrl-2009}:
\begin{align}
\sum_{k=0}^{\infty}\frac{4k+1}{81^k}\gamma_k=\frac{3\sqrt{3}}{2\pi}.\label{aa-1}
\end{align}

The above interesting example motivates us to prove the following supercongruence, which was originally conjectured by Zudilin \cite[(33)]{zudilin-jnt-2009}.
\begin{thm}\label{t-1}
For any prime $p\ge 5$, we have
\begin{align}
\sum_{k=0}^{p-1}\frac{4k+1}{81^k}\gamma_k \equiv \left(\frac{-3}{p}\right)p\pmod{p^3},\label{aa-2}
\end{align}
where $\left(\frac{\cdot}{p}\right)$ denotes the Legendre symbol.
\end{thm}

 The supercongruence \eqref{aa-2} may be regarded as a $p$-adic analogue of \eqref{aa-1}. In the past two decades, Ramanujan-type series for $1/\pi$ as well as related supercongruences and $q$-supercongruences have attracted many experts' attention (see, for instance, \cite{ccl-am-2004,cv-mrl-2009,cwz-ijm-2013,guo-aam-2020,gl-itsf-2020,gs-rm-2020,
gz-am-2019,liu-jsc-2019,liu-jmaa-2020,lp-a-2020,sunzw-scm-2011,sunzw-era-2020,van-b-1997,zudilin-jnt-2009}).

The second result of this paper consists of the following two related supercongruences involving
the numbers $\gamma_n$, which were originally conjectured by Z.-H. Sun \cite[Conjecture 6.8]{sunzh-itsf-2015}.
\begin{thm}\label{t-2}
For any prime $p\ge 5$, we have
\begin{align}
\sum_{k=0}^{p-1}(4k+3)\gamma_k \equiv 3\left(\frac{-3}{p}\right)p\pmod{p^3}.\label{aa-3}
\end{align}
\end{thm}

\begin{thm}\label{t-3}
For any prime $p\ge 5$, we have
\begin{align}
\sum_{k=0}^{p-1}\frac{2k+1}{(-9)^k}\gamma_k\equiv \left(\frac{-3}{p}\right)p\pmod{p^3}.\label{aa-4}
\end{align}
\end{thm}

We remark that congruence properties for the Almkvist--Zudilin numbers have been widely investigated by Amdeberhan and Tauraso \cite{at-aa-2016}, Chan, Cooper and Sica \cite{ccs-ijnt-2010}, and Z.-H. Sun \cite{sunzh-itsf-2015,sunzh-a-20-2,sunzh-a-20-4}.

The rest of the paper is organized as follows. In Section 2, we recall some necessary combinatorial identities involving harmonic numbers and prove a preliminary congruence. The proofs of Theorems \ref{t-1}--\ref{t-3} are presented in Sections 3--5, respectively.

\section{Preliminary results}
Let
\begin{align*}
H_n=\sum_{j=1}^n\frac{1}{j}
\end{align*}
denote the $n$th harmonic number. The Fermat quotient of an integer $a$ with respect to an odd prime $p$ is given by $q_p(a)=(a^{p-1}-1)/p$.

In order to prove Theorems \ref{t-1} and \ref{t-2}, we need the following two lemmas.
\begin{lem}
For any non-negative integer $n$, we have
\begin{align}
&\sum_{i=0}^{n}(-1)^i{n\choose i}{n+i\choose i}=(-1)^n,\label{bb-1}\\
&\sum_{i=0}^{n}(-1)^i{n\choose i}{n+i\choose i}H_i=2(-1)^nH_n,\label{bb-2}\\
&\sum_{i=0}^{n}(-1)^i{n\choose i}{n+i\choose i}H_{n+i}=2(-1)^nH_n.\label{bb-3}
\end{align}
\end{lem}
In fact, such identities can be discovered and proved by the symbolic summation package {\tt Sigma} developed by Schneider \cite{schneider-slc-2007}. One can also refer to \cite{liu-jsc-2019} for the same approach to finding and proving identities of this type. For human proofs of \eqref{bb-1}--\eqref{bb-3}, one refers to \cite{prodinger-integers-2008}.

\begin{lem}
For any prime $p\ge 5$, we have
\begin{align}
\sum_{k=0}^{p-1}
\frac{(3k)!}{3^{3k} k!^3}\left(H_{3k}-H_k\right)\equiv \left(\frac{-3}{p}\right)q_p(3)\pmod{p}.
\label{bb-4}
\end{align}
\end{lem}

{\noindent\it Proof.}
Note that
\begin{align}
\sum_{k=0}^{p-1}
\frac{(3k)!}{3^{3k} k!^3}\left(3H_{3k}-H_k\right)
=\sum_{k=0}^{p-1}\frac{(1/3)_k(2/3)_k}{(1)_k^2}\sum_{j=0}^{k-1}
\left(\frac{1}{1/3+j}+\frac{1}{2/3+j}\right).
\label{bb-5}
\end{align}
Recall the following identity due to Tauraso \cite[Theorem 1]{tauraso-integers-2012}:
\begin{align}
\frac{(1/3)_k(2/3)_k}{(1)_k^2}\sum_{j=0}^{k-1}\left(\frac{1}{1/3+j}+\frac{1}{2/3+j}\right)
=\sum_{j=0}^{k-1}\frac{(1/3)_j(2/3)_j}{(1)_j^2}\cdot \frac{1}{k-j}.
\label{bb-6}
\end{align}
Substituting \eqref{bb-6} into \eqref{bb-5} and exchanging the summation order gives
\begin{align}
\sum_{k=0}^{p-1}
\frac{(3k)!}{3^{3k} k!^3}\left(3H_{3k}-H_k\right)
&=\sum_{k=0}^{p-1}\sum_{j=0}^{k-1}\frac{(1/3)_j(2/3)_j}{(1)_j^2}\cdot \frac{1}{k-j}\notag\\[10pt]
&=\sum_{j=0}^{p-2}\frac{(1/3)_j(2/3)_j}{(1)_j^2}H_{p-1-j}\notag\\[10pt]
&\equiv \sum_{j=0}^{p-2}\frac{(3j)!}{3^{3j} j!^3}H_{j}\pmod{p},\label{bb-7}
\end{align}
where we have utilized the fact that $H_{p-1-j}\equiv H_j\pmod{p}$.
By \eqref{bb-7}, we obtain
\begin{align}
\sum_{k=0}^{p-1}\frac{(3k)!}{3^{3k} k!^3}\left(H_{3k}-H_k\right)
&\equiv \frac{1}{3}\left(\sum_{k=0}^{p-2}\frac{(3k)!}{3^{3k} k!^3}H_{k}-2\sum_{k=0}^{p-1}\frac{(3k)!}{3^{3k} k!^3}H_k\right)\notag\\[10pt]
&\equiv -\frac{1}{3}\sum_{k=0}^{\lfloor p/3\rfloor}\frac{(3k)!}{3^{3k} k!^3}H_{k} \pmod{p},
\label{bb-8}
\end{align}
because $(3k)!\equiv 0\pmod{p}$ for $k>\lfloor p/3\rfloor$.

Let $m=\lfloor p/3\rfloor$. From \cite[Lemma 2.3]{at-aa-2016}, we see that for $0\le k\le m$,
\begin{align}
\frac{(3k)!}{3^{3k} k!^3}\equiv (-1)^k{m\choose k}{m+k\choose k}\pmod{p}.
\label{bb-9}
\end{align}
It follows from \eqref{bb-2}, \eqref{bb-8} and \eqref{bb-9} that
\begin{align*}
\sum_{k=0}^{p-1}\frac{(3k)!}{3^{3k} k!^3}\left(H_{3k}-H_k\right)
&\equiv -\frac{1}{3}\sum_{k=0}^m(-1)^k{m\choose k}{m+k\choose k}H_k\pmod{p}\\[10pt]
&=-\frac{2(-1)^m}{3}H_m.
\end{align*}

Finally, noting
\begin{align}
(-1)^{\lfloor p/3\rfloor}=\left(\frac{-3}{p}\right),\label{bb-10}
\end{align}
and the following congruence \cite[page 359]{lehmer-am-1938}:
\begin{align}
H_{\lfloor p/3\rfloor}\equiv -\frac{3}{2}q_p(3) \pmod{p^2},\label{bb-11}
\end{align}
we complete the proof of \eqref{bb-4}.
\qed

\section{Proof of Theorem \ref{t-1}}
We begin with the transformation formula due to Chan and Zudilin \cite[Corollary 4.3]{cz-m-2010}:
\begin{align}
\gamma_n=\sum_{i=0}^n{2i\choose i}^2{4i\choose 2i}{n+3i\choose 4i}(-3)^{3(n-i)}.\label{cc-1}
\end{align}
Using \eqref{cc-1} and exchanging the summation order, we obtain
\begin{align}
\sum_{k=0}^{p-1}\frac{4k+1}{81^k}\gamma_k
&=\sum_{k=0}^{p-1}\frac{4k+1}{81^k}\sum_{i=0}^k{2i\choose i}^2{4i\choose 2i}{k+3i\choose 4i}(-3)^{3(k-i)}\notag\\[10pt]
&=\sum_{i=0}^{p-1}\frac{1}{(-3)^{3i}}{2i\choose i}^2{4i\choose 2i}\sum_{k=i}^{p-1}\frac{4k+1}{(-3)^k}{k+3i\choose 4i}.\label{cc-2}
\end{align}
Note that
\begin{align}
\sum_{k=i}^{n-1}\frac{4k+1}{(-3)^k}{k+3i\choose 4i}=(n-i){n+3i\choose 4i}(-3)^{1-n},\label{cc-3}
\end{align}
which can be easily proved by induction on $n$. Combining \eqref{cc-2} and \eqref{cc-3} gives
\begin{align}
\sum_{k=0}^{p-1}\frac{4k+1}{81^k}\gamma_k
=3^{1-p}\sum_{i=0}^{p-1}\frac{p-i}{(-3)^{3i}}{2i\choose i}^2{4i\choose 2i}{p+3i\choose 4i}.\label{cc-4}
\end{align}

Furthermore, we have
\begin{align*}
&(-1)^i(p-i){2i\choose i}^2{4i\choose 2i}{p+3i\choose 4i}\\[10pt]
&=\frac{(-1)^ip(p+3i)\cdots (p+1)(p-1)\cdots (p-i)}{i!^4}\\[10pt]
&\equiv \frac{p(3i)!}{i!^3}\left(1+p\left(H_{3i}-H_i\right)\right)\pmod{p^3}.
\end{align*}
Thus,
\begin{align*}
\sum_{k=0}^{p-1}\frac{4k+1}{81^k}\gamma_k
&\equiv 3^{1-p}p\sum_{i=0}^{p-1}
\frac{(3i)!}{3^{3i} i!^3}\left(1+p\left(H_{3i}-H_i\right)\right)\pmod{p^3}.
\end{align*}

Finally, noting \eqref{bb-4} and Mortenson's supercongruence \cite[(1.2)]{mortenson-tams-2003}:
\begin{align*}
\sum_{i=0}^{p-1}
\frac{(3i)!}{3^{3i} i!^3}
\equiv \left(\frac{-3}{p}\right) \pmod{p^2},
\end{align*}
we arrive at
\begin{align*}
\sum_{k=0}^{p-1}\frac{4k+1}{81^k}\gamma_k
&\equiv p\left(\frac{-3}{p}\right)\left(3^{1-p}+3^{1-p}pq_p(3)\right)\pmod{p^3}\\[10pt]
&=p\left(\frac{-3}{p}\right),
\end{align*}
as desired.

\section{Proof of Theorem \ref{t-2}}
Recall the following transformation formula \cite[(5.1)]{sunzh-itsf-2015}:
\begin{align}
\gamma_n=\sum_{i=0}^{\lfloor n/3 \rfloor}{2i\choose i}^2{4i\choose 2i}{n+i\choose 4i}(-3)^{n-3i}.\label{dd-1}
\end{align}
By \eqref{dd-1}, we have
\begin{align}
\sum_{k=0}^{p-1}(4k+3)\gamma_k
&=\sum_{k=0}^{p-1}(4k+3)\sum_{i=0}^{\lfloor k/3 \rfloor}{2i\choose i}^2{4i\choose 2i}{k+i\choose 4i}(-3)^{k-3i}\notag\\[10pt]
&=\sum_{i=0}^{p-1}\frac{1}{(-3)^{3i}}{2i\choose i}^2{4i\choose 2i}\sum_{k=i}^{p-1}(-3)^k(4k+3){k+i\choose 4i}.\label{dd-2}
\end{align}
It can be easily proved by induction on $n$ that
\begin{align}
\sum_{k=i}^{n-1}(-3)^k(4k+3){k+i\choose 4i}=3(n-3i){n+i\choose 4i}(-3)^{n-1}.\label{dd-3}
\end{align}
It follows from \eqref{dd-2} and \eqref{dd-3} that
\begin{align*}
\sum_{k=0}^{p-1}(4k+3)\gamma_k
=3^{p}\sum_{i=0}^{\lfloor p/3\rfloor}\frac{p-3i}{(-3)^{3i}}{2i\choose i}^2{4i\choose 2i}{p+i\choose 4i}.
\end{align*}

Note that
\begin{align*}
&(-1)^i(p-3i){2i\choose i}^2{4i\choose 2i}{p+i\choose 4i}\\[10pt]
&=\frac{(-1)^ip(p+i)\cdots(p+1)(p-1)\cdots(p-3i)}{i!^4}\\[10pt]
&\equiv \frac{p(3i)!}{i!^3}\left(1-p\left(H_{3i}-H_i\right)\right)\pmod{p^3}.
\end{align*}
Thus,
\begin{align}
\sum_{k=0}^{p-1}(4k+3)\gamma_k
&\equiv 3^{p}p
\sum_{i=0}^{\lfloor p/3\rfloor}\frac{(3i)!}{3^{3i} i!^3}\left(1-p\left(H_{3i}-H_i\right)\right)\pmod{p^3}.\label{dd-4}
\end{align}

Let $m=\lfloor p/3\rfloor$. Since
\begin{align*}
\frac{(3i)!}{3^{3i} i!^3}=(-1)^i{-1/3\choose i}{-1/3+i\choose i}
=(-1)^i{-2/3\choose i}{-2/3+i\choose i},
\end{align*}
we have
\begin{align}
\frac{(3i)!}{3^{3i} i!^3}
&=(-1)^i{m-p/3\choose i}{m-p/3+i\choose i}\notag\\[10pt]
&=\frac{(-1)^i(m+i-p/3)\cdots (m-i+1-p/3)}{i!^2}\notag\\[10pt]
&\equiv (-1)^i{m\choose i}{m+i\choose i}\left(1-\frac{p}{3}\left(H_{m+i}-H_{m-i}\right)\right)\pmod{p^2}.
\label{dd-5}
\end{align}
Substituting \eqref{dd-5} into the right-hand side of \eqref{dd-4} gives
\begin{align}
&\sum_{k=0}^{p-1}(4k+3)\gamma_k\notag\\
&\equiv 3^{p}p
\sum_{i=0}^{m}(-1)^i{m\choose i}{m+i\choose i}\notag\\
&\times \left(1-\frac{p}{3}\left(H_{m+i}-H_{m-i}+3H_{3i}-3H_i\right)\right)
\pmod{p^3}.\label{dd-6}
\end{align}

Furthermore, we have
\begin{align}
H_{3i}&=\frac{1}{3}\left(H_i+\sum_{j=1}^i\frac{1}{j-1/3}+\sum_{j=1}^i\frac{1}{j-2/3}\right)\notag\\[10pt]
&\equiv \frac{1}{3}\left(H_i+\sum_{j=1}^i\frac{1}{m+j}-\sum_{j=1}^i\frac{1}{m+1-j}\right)\pmod{p}\notag\\[10pt]
&=\frac{1}{3}\left(H_i+H_{m+i}+H_{m-i}-2H_m\right).\label{dd-7}
\end{align}
It follows from \eqref{bb-1}--\eqref{bb-3}, \eqref{dd-6} and \eqref{dd-7} that
\begin{align*}
&\sum_{k=0}^{p-1}(4k+3)\gamma_k\\
&\equiv 3^{p}p
\sum_{i=0}^{m}(-1)^i{m\choose i}{m+i\choose i}\left(1-\frac{2p}{3}\left(H_{m+i}-H_m-H_i\right)\right)
\pmod{p^3}\\
&=3^{p}p(-1)^m\left(1+\frac{2p}{3}H_m\right).
\end{align*}

Finally, using \eqref{bb-10} and \eqref{bb-11}, we obtain
\begin{align*}
\sum_{k=0}^{p-1}(4k+3)\gamma_k
&\equiv 3p\left(\frac{-3}{p}\right)3^{p-1}\left(2-3^{p-1}\right)\\[10pt]
&=3p\left(\frac{-3}{p}\right)\left(1-\left(3^{p-1}-1\right)^2\right)\\[10pt]
&\equiv 3p\left(\frac{-3}{p}\right)\pmod{p^3},
\end{align*}
where we have used the Fermat's little theorem in the last step.

\section{Proof of Theorem \ref{t-3}}
Recall the following transformation formula \cite[Lemma 4.1]{sunzh-itsf-2015}:
\begin{align}
\gamma_n=\sum_{i=0}^{n}(-9)^{n-i}{2i\choose i}{n+i\choose 2i}\sum_{j=0}^i{i\choose j}^2{2j\choose j}.\label{ee-1}
\end{align}
Let
\begin{align*}
g_n=\sum_{k=0}^n{n\choose k}^2{2k\choose k}.
\end{align*}
By \eqref{ee-1}, we have
\begin{align}
\sum_{k=0}^{p-1}\frac{2k+1}{(-9)^k}\gamma_k
&=\sum_{k=0}^{p-1}\frac{2k+1}{(-9)^k}\sum_{i=0}^{k}(-9)^{k-i}{2i\choose i}{k+i\choose 2i}g_i\notag\\[10pt]
&=\sum_{i=0}^{p-1}\frac{g_i}{(-9)^i}\sum_{k=i}^{p-1}(2k+1){k+i\choose 2i}{2i\choose i}.\label{ee-2}
\end{align}
Note that
\begin{align}
\sum_{k=i}^{n-1}(2k+1){k+i\choose 2i}{2i \choose i}
=\frac{n^2}{i+1}{n-1\choose i}{n+i\choose i},\label{ee-3}
\end{align}
which can be proved by induction on $n$. It follows from \eqref{ee-2} and \eqref{ee-3} that
\begin{align*}
\sum_{k=0}^{p-1}\frac{2k+1}{(-9)^k}\gamma_k
=p^2\sum_{i=0}^{p-1}\frac{g_i}{(-9)^i(i+1)} {p-1\choose i}{p+i\choose i}.
\end{align*}
Since
\begin{align*}
{p-1\choose i}{p+i\choose i}\equiv (-1)^i\pmod{p^2},
\end{align*}
we have
\begin{align}
\sum_{k=0}^{p-1}\frac{2k+1}{(-9)^k}\gamma_k\equiv p^2\sum_{i=0}^{p-1}\frac{g_i}{9^i(i+1)}\pmod{p^3}.
\label{ee-4}
\end{align}

From \cite[Lemma 2.7]{jv-rj-2010}, we see that for $0\le i\le p-1$,
\begin{align*}
\frac{g_i}{9^i}\equiv \left(\frac{-3}{p}\right)g_{p-1-i}\pmod{p},
\end{align*}
and so
\begin{align*}
\sum_{i=0}^{p-2}\frac{g_i}{9^i(i+1)}
&\equiv \left(\frac{-3}{p}\right)\sum_{i=0}^{p-2}\frac{g_{p-1-i}}{i+1}\\[10pt]
&=\left(\frac{-3}{p}\right)\sum_{i=1}^{p-1}\frac{g_{i}}{p-i}\\[10pt]
&\equiv -\left(\frac{-3}{p}\right)\sum_{i=1}^{p-1}\frac{g_{i}}{i}\pmod{p}.
\end{align*}
Using the congruence \cite[(1.8)]{sunzw-rj-2016}:
\begin{align*}
\sum_{i=1}^{p-1}\frac{g_{i}}{i}\equiv 0\pmod{p},
\end{align*}
we obtain
\begin{align}
\sum_{i=0}^{p-2}\frac{g_i}{9^i(i+1)}\equiv 0\pmod{p}.\label{ee-5}
\end{align}
Furthermore, combining \eqref{ee-4} and \eqref{ee-5} gives
\begin{align*}
\sum_{k=0}^{p-1}\frac{2k+1}{(-9)^k}\gamma_k\equiv \frac{pg_{p-1}}{9^{p-1}}\pmod{p^3}.
\end{align*}

By \cite[Lemma 3.2]{sunzw-rj-2016}, we have
\begin{align*}
g_{p-1}\equiv \left(\frac{-3}{p}\right)\left(2\cdot3^{p-1}-1\right)\pmod{p^2},
\end{align*}
and so
\begin{align*}
\sum_{k=0}^{p-1}\frac{2k+1}{(-9)^k}\gamma_k
&\equiv p\left(\frac{-3}{p}\right)
\left(1-\frac{\left(3^{p-1}-1\right)^2}{9^{p-1}}\right)\\[10pt]
&\equiv p\left(\frac{-3}{p}\right)\pmod{p^3},
\end{align*}
where we have utilized the Fermat's little theorem.

{\noindent\bf Remark.}
Z.-H. Sun \cite[Conjecture 6.8]{sunzh-itsf-2015} also conjectured a companion supercongruence of
\eqref{aa-4}:
\begin{align}
\sum_{k=0}^{p-1}\frac{2k+1}{9^k}\gamma_k \equiv \left(\frac{-3}{p}\right)p\pmod{p^3}.
\label{ee-6}
\end{align}
In a similar way, by using \eqref{ee-1} and the following identity:
\begin{align*}
\sum_{k=i}^{n-1}(-1)^k(2k+1){k+i\choose 2i}{2i \choose i}
=(-1)^{n-1}n{n-1\choose i}{n+i\choose i},
\end{align*}
we can show that
\begin{align*}
\sum_{k=0}^{p-1}\frac{2k+1}{9^k}\gamma_k
\equiv p\sum_{i=0}^{p-1}\frac{g_i}{9^i}\pmod{p^3}.
\end{align*}
Thus, the conjectural supercongruence \eqref{ee-6} is equivalent to
\begin{align*}
\sum_{i=0}^{p-1}\frac{g_i}{9^i}\equiv \left(\frac{-3}{p}\right)\pmod{p^2},
\end{align*}
which was originally conjectured by Z.-W. Sun \cite[Remark 1.1]{sunzw-rj-2016}.

\vskip 5mm \noindent{\bf Acknowledgments.}
This work was supported by the National Natural Science Foundation of China (grant 11801417).

\end{document}